\documentclass[reqno, twoside,a4paper]{amsart}
\usepackage{style}

\usepackage{graphicx} 

\begin{document}

\title[Orthogonal polynomials]{On rate of convergence for universality limits}
\author{Roman Bessonov}

\address{
\begin{flushleft}
bessonov@pdmi.ras.ru\\\vspace{0.1cm}
St.\,Petersburg State University\\  
Universitetskaya nab. 7-9, 199034 St.\,Petersburg, RUSSIA\\
\vspace{0.1cm}
St.\,Petersburg Department of Steklov Mathematical Institute\\ Russian Academy of Sciences\\
Fontanka 27, 191023 St.Petersburg,  RUSSIA\\
\end{flushleft}
}
\thanks{
The work is supported by grant RScF 19-71-30002-$\Pi$ of the Russian Science Foundation.
}


\begin{abstract} 
Given a probability measure $\mu$ on the unit circle $\T$, consider the reproducing kernel $k_{\mu,n}(z_1, z_2)$ in the space of polynomials of degree at most $n-1$ with the $L^2(\mu)$--inner product. Let $u, v \in \C$. It is known that under  mild assumptions on $\mu$ near $\zeta \in \mathbb{T}$, the ratio $k_{\mu,n}(\zeta e^{u/n}, \zeta e^{v/n})/k_{\mu,n}(\zeta, \zeta)$ converges to a universal limit $S(u, v)$ as $n \to \infty$. We give an estimate for the rate of this convergence for measures $\mu$ with finite logarithmic integral. 
\end{abstract}

\maketitle

\section{Introduction}
Consider a probability measure $\mu$ on the unit circle $\T = \{z \in \C: \; |z| = 1\}$ of the complex plane, $\C$. Assume that the support of $\mu$ is an infinite subset of $\T$, so that monomials $z^k$, $k \ge 0$, are linearly independent in $L^2(\mu)$. For each integer $n \ge 1$, the set of polynomials of degree at most $n-1$,
$$
\mathcal P_{n} = \spn\{z^k,\; k=0,\ldots, n-1\},
$$
can be viewed as the $n$-dimensional Hilbert space of analytic functions with respect to  $L^2(\mu)$-inner product. Denote by $k_{\mu,n}(z_1, z_2)$ the reproducing kernel at a point $z_2 \in \C$ in this space, i.e., $k_{\mu,n}(\cdot, z_2) \in \mathcal P_{n}$ and
$$
(p, k_{\mu,n}(\cdot, z_2))_{L^2(\mu)} = p(z_2), \qquad p \in \mathcal P_{n}.
$$  
If $\mu = m$ is the Lebesgue measure on $\T$ normalized by $m(\T) = 1$, the reproducing kernel has the following form:
$$
k_{m,n}(z_1, z_2) = \frac{1 - \ov{z_2}^{n} z_1^{n}}{1 - \ov{z_2} z_1}.
$$
One might check that if $z_1 = \zeta e^{u/n}$, $z_2 = \zeta e^{v/n}$ for some $\zeta \in \T$ and $u, v \in \C$ (equivalently, $z_1$, $z_2$ are at a distance comparable to $1/n$ from $\zeta$), then
$$
\frac{k_{m,n}(z_1, z_2)}{k_{m,n}(\zeta, \zeta)} = \frac{1 - \ov{z_2}^n z_1^n}{n(1 - \ov{z_2} z_1)}  = 
\frac{e^{u+\bar v} - 1}{u+\bar v} + O\left(\frac{1}{n}\right),
$$
where the remainder is uniform in $(u,v)$  on compact subsets of $\C \times \C$. Such kind of behaviour of reproducing kernels is universal:
under mild assumptions on a measure $\mu$ near $\zeta \in \T$, we have 
$$
\frac{k_{\mu,n}(\zeta e^{u/n},\, \zeta e^{v/n})}{k_{\mu,n}(\zeta, \zeta)} \to \frac{e^{u+\bar v} - 1}{u+\bar v}, \qquad n \to +\infty,
$$ 
uniformly in $(u, v)$ on compact subsets of $\C \times \C$. Universality of the limiting behaviour of ratios of reproducing kernels attracted major attention in recent years. Several essentially different approaches were developed. Let us mention some of them. First papers dealt with real analytic weights and used the Riemann-Hilbert method, see, e.g., P.~Deift~\cite{Deift} or A.~Kuijlaars, M.~Vanlessen \cite{KV02}. D.~Lubinsky \cite{Lubinsky} found a way of reducing a wide class of universality problems to the study of asymptotic behaviour of $k_{\mu,n}(z,z)$, $z \in \T$. The latter asymptotic behaviour has been previously identified for general measures of Szeg\H{o} class by A.~M\'at\'e, P.~Nevai and V.~Totik in \cite{MNT91}. Global Szeg\H{o} condition has been weakened to the local one by E.~Findley in \cite{Findley}. Another approach, also pioneered by D.~Lubinsky \cite{Lubinsky08}, is based on compactness of normal families of entire functions and properties of $\frac{\sin x}{x}$ kernel. An overview of this approach and further results can be found in B.~Simon \cite{Simon08} and V.~Totik \cite{Totik}. Recently, B.~Eichinger, M.~Lukic, and B.~Simanek \cite{ELS} found yet another approach to universality based on spectral theory of canonical Hamiltonian systems. While this approach gives extremely general results (even the local Szeg\H{o} condition can be omitted), it also involves a compactness argument as an essential element of the proof. Most of mentioned papers deal with measures on subsets of the real line due to motivation in the theory of random matrices. However, even in the simplified setting of measures on the unit circle, estimates for the rate of convergence 
\begin{equation}\label{eq27}
\frac{k_{\mu,n}(z_1, z_2)}{k_{\mu,n}(\zeta, \zeta)} \to \frac{e^{u+\bar v} - 1}{u+\bar v}   
\end{equation}
are missing in the literature. In fact, the rate of convergence in \eqref{eq27} is not known even in the case where $\mu$ is an absolutely continuous measure on $\T$ with a smooth non-vanishing weight $w$. Indeed, compactness arguments, widely used for proving universality, cannot give bounds for the rate of convergence. 

\medskip

As an additional motivation of this work, we mention that Poltoratski \cite{Polt21} recently used universality in the proof of convergence of certain nonlinear Fourier transform (NLFT), and a subsequent development of this area, e.g., bounds for NLFT maximal operators, will require estimates for the convergence of universality limits.

\medskip

In this paper, we estimate the rate of convergence in \eqref{eq27} for probability measures on the unit circle with finite logarithmic integral. For this we use an entropy function of a measure -- a powerful instrument that recently found several applications in inverse problems \cite{BD2017}, \cite{BD2019}, scattering theory \cite{BD2022}, and orthogonal polynomials \cite{BD21}, \cite{B21}. To be precise, let $\mu$ be a probability measure on $\T$, and let $\mu = w\,dm + \mus$ be its Radon-Nikodym decomposition into the absolutely continuous and singular parts. The measure $\mu$ is said to belong to the Szeg\H{o} class $\szc$ if its logarithmic integral is finite:
$$
\int_{\T}\log w\, dm > -\infty.
$$
Since $\log x \le x$ for all $x>0$, the latter condition is equivalent to $\log w \in L^1(m)$. Let $\D = \{z \in \C:\; |z| < 1\}$ denote the open unit disk. For a measure $\mu \in \szc$, its entropy function is given by
$$
\K_{\mu}(z) = \log\left(\int_{\T}\frac{1-|z|^2}{|1-\bar \xi z|^2}\,d\mu(\xi)\right) - \int_{\T}\log w(\xi)\frac{1-|z|^2}{|1-\bar \xi z|^2}\,dm(\xi), \qquad z \in \D.
$$
The function $\K_{\mu}$ is nonnegative in $\D$ by Jensen's inequality. Moreover, for $m$-almost all $\zeta \in \T$, we have $\K_{\mu}(z) \to 0$ as $z$ non-tangentially approaches $\zeta$. This follows from well-known properties of the Poisson kernel: we have
\begin{equation}\label{eq22}
\lim_{r \to 1}\int_{\T}\frac{1-r^2}{|1-\bar \xi r\zeta|^2}\,d\mus(\xi) = 0
\end{equation}
for $m$-almost all $\zeta \in \T$, and 
\begin{align}
&\lim_{r \to 1}\int_{\T}w(\xi)\frac{1-r^2}{|1-\bar \xi r\zeta|^2}\,dm(\xi) \to w(\zeta), \label{eq23}\\
&\lim_{r \to 1}\int_{\T}\log w(\xi)\frac{1-r^2}{|1-\bar \xi r\zeta|^2}\,dm(\xi) \to \log w(\zeta), \label{eq24}
\end{align}
at each Lebesgue point $\zeta$ of functions $w, \log w \in L^1(m)$. In case \eqref{eq22}-\eqref{eq24} are satisfied, we have 
$\K_\mu(r\zeta) \to 0$ as $r \to 1$.  Let
$B(\zeta, r) = \{z \in \C:\; |z-\zeta| \le r\}$. 
The following theorem is the main result of the paper.
\begin{Thm}\label{t1}
Let $\mu \in \szc$, $A \ge 1$, $n \ge 10A$, $\zeta \in \T$. There exists $\eps_0 > 0$ depending only on $A$, such that if $z_{1,2} \in B(\zeta, A/n)$ and $\K_{\mu}(\rho\zeta) \le \eps_0$ for all $\rho \in [1-A/n, 1)$, then 
\begin{equation}\label{eq11}
\left|\frac{k_{\mu, n}(z_1, z_2)}{k_{\mu, n}(\zeta, \zeta)} - \frac{1 - \ov{z_2}^n z_1^n}{n(1 - \ov{z_2} z_1)}\right| \le c e^{4A}\sup_{\rho \in [1-\delta, 1) }\sqrt{\K_{\mu}(\rho\zeta)}, 
\end{equation}
where $\delta = \max\limits_{k=1,2}|z_{k} - \zeta|$, and the constant $c>0$ is absolute.
\end{Thm}
Note that $\rho \in [1- A/n, 1)$ in Theorem \ref{t1} tends to $1$ as $n \to \infty$, therefore, the right hand side of \eqref{eq11} tends to zero for all $\zeta$ satisfying \eqref{eq22}-\eqref{eq24}. This gives a nontrivial bound for the rate of convergence in \eqref{eq27} for Lebesgue almost all $\zeta \in \T$. If $\mu$ has some regularity in a neighbourhood of $\zeta \in \T$, its entropy function can be explicitly estimated. For functions $f, g > 0$, we use notation $f \lesssim g$ (resp., $f \gtrsim g$) if $f \le cg$ (resp., $f \le cg$) for some constant $c$, and $f \sim g$ if both relations $f \lesssim g$ and $f \gtrsim g$ are satisfied.  
\begin{Thm}\label{t2}
Let $\mu = w\,dm$ be an absolutely continuous probability measure in $\szc$ such that $w$ is positive and continuous in a neighbourhood $I \subset \T$ of  $\zeta \in \T$. Assume, moreover, that $|w(\xi) -  w(\zeta)| \sim |\xi - \zeta|^{s}$ for all $\xi \in I$ and some $s> 0$. Then we have 
$$
\K_{\mu}(\rho\zeta) \sim 
\begin{cases}
1-\rho, & s \in (1/2, +\infty),\\
(1-\rho)|\log(1-\rho)|,  & s =1/2,\\
(1-\rho)^{2s}, & s \in (0, 1/2),
\end{cases}
$$ 
for $\rho \in (0, 1)$ close enough to $1$. The constants involved depend on $s$, the diame\-ter of~$I$, the value $w(\zeta)$, and the constants in the relation $|w(\xi) -  w(\zeta)| \sim |\xi - \zeta|^{s}$.
\end{Thm}
\noindent Let $\lambda \in \D$. For the absolutely continuous probability measure $\mu = \frac{1 - |\lambda|^2}{|1-\lambda \xi|^2}\,dm$, we have 
\begin{align*}
\K(\mu, z) 
&= \log\Re\left(\frac{1+\lambda z}{1-\lambda z}\right) - \log\left(\frac{1-|\lambda|^2}{|1-\lambda z|^2}\right)\\
&= \log\frac{1-|\lambda z|^2}{|1-\lambda z|^2} - \log\left(\frac{1-|\lambda|^2}{|1-\lambda z|^2}\right) = \log\frac{1-|\lambda z|^2}{1-|\lambda|^2}, 
\end{align*}
due to the fact that integration against the Poisson kernel corresponds to harmonic continuation into the unit disk. We see that $\K(\mu, (1-1/n)\zeta) \sim 1/n$ as $n \to \infty$. Note that this agrees with bounds in Theorem \ref{t2} (we have $s=1$ for this measure). By Theorem \ref{t1}, we then have
$$
\left|\frac{k_{\mu, n}(z_1, z_2)}{k_{\mu, n}(\zeta, \zeta)} - \frac{1 - \ov{z_2}^n z_1^n}{n(1 - \ov{z_2} z_1)}\right| \lesssim \frac{1}{\sqrt{n}}
$$
for all $z_1$, $z_2$ in $B(\zeta, 1/n)$ and large enough $n \ge 0$, uniformly in $\zeta \in \T$. As we will see in Section \ref{s4}, in fact 
$$
\sup_{|z_{1,2}-\zeta| \le 1/n}\left|\frac{k_{\mu, n}(z_1, z_2)}{k_{\mu, n}(\zeta, \zeta)} - \frac{1 - \ov{z_2}^n z_1^n}{n(1 - \ov{z_2} z_1)}\right| \sim \frac{1}{n}.
$$
This shows that the bound in Theorem \ref{t1} is not sharp for smooth measures. It seems, however, that this bound cannot be improved in the setting of the whole class $\szc$ of measures with finite logarithmic integral, i.e., there is a measure $\mu \in \szc$ such that 
$$
\sup_{|z_{1,2} - 1| \le 1/n}\left|\frac{k_{\mu, n}(z_1, z_2)}{k_{\mu, n}(1, 1)} - \frac{1 - \ov{z_2}^n z_1^n}{n(1 - \ov{z_2} z_1)}\right| \gtrsim  \sqrt{\K_{\mu}(1-1/n)}, 
$$
for all $n \ge 0$. In Section \ref{s5}, we consider the absolutely continuous measures $w_s\,dm$ such that $w_s(e^{i\theta}) = c_{s}e^{|\theta|^{s}}$, $\theta \in [-\pi, \pi]$,  where the constant $c_s$ is chosen so that $\int_{\T} w_s \,dm = 1$. By Theorem \ref{t2}, we have
$$
\K_{w_s}(1 - 1/n) \sim n^{-2s}, \qquad s \in (0, 1/2).
$$
We demonstrate numerically that for $x_n = 1- n^{-1}$, and fixed $s \in (0, 1/2)$, we have 
$$
\left|\frac{k_{w_s, n}(x_n, x_n)}{k_{w_s, n}(1, 1)} - \frac{1 - |x_n|^{2n}}{n(1 - |x_n|^2)}\right| \gtrsim  n^{-s}.
$$
In other words, for each $s \in (0,1/2)$, we have
$$
\sup_{|z_{1,2} - 1| \le 1/n}\left|\frac{k_{w_s, n}(z_1, z_2)}{k_{w_s, n}(1, 1)} - \frac{1 - \ov{z_2}^n z_1^n}{n(1 - \ov{z_2} z_1)}\right| \gtrsim  \sqrt{\K_{w_s}(1-1/n)}, 
$$
and estimate \eqref{eq11} in Theorem \ref{t1} is sharp on this class of examples. It remains an open problem to give a mathematical proof of this fact. 
\section{Proof of Theorem \ref{t1}}
Let $\mu$ be a probability measure supported on an infinite subset of $\T$, and let $\{\phi_n\}_{n \ge 0}$ be the family of its orthonormal polynomials obtained by Gram-Schmidt orthogonalization of monomials $z^n$, $n \ge 0$, in $L^2(\mu)$. 
For a polynomial $p$ of degree $n$, we set $p^*(z) = z^n\ov{p(1/
\bar z)}$. Note that $p$ is also a polynomial of degree at most $n$. The polynomials $\{\phi_n^*\}_{n \ge 0}$ are called reflected orthonormal polynomials. We have the following recurrence relation (see formula $(1.5.25)$, page 58, in \cite{Simonbook}): 
$$
\phi_{n} = \frac{z\phi_{n-1} - 
\ov{a_{n-1}}\phi^{*}_{n-1}}{\sqrt{1-|a_{n-1}|^2}}, \qquad n \ge 1. 
$$
Here the recurrence coefficients, $a_n$, $n\ge 0$, belong to $\D$. It is also known (see Theorem 2.2.7, page 124, in \cite{Simonbook}) that the reproducing kernel in the $n$-dimensional space Hilbert space $(\mathcal P_n, (\cdot, \cdot)_{L^2(\mu)})$ at $z_2 \in \C$ is given by
\begin{equation}\label{eq12}
k_{\mu, n}(z_1, z_2) =
\sum_{k=0}^{n-1}\ov{\phi_{k}(z_2)}\phi_{k}(z_1)
= \frac{\ov{\phi_{n}^*(z_2)}\phi_{n}^*(z_1)
-
\ov{\phi_{n}(z_2)}\phi_{n}(z_1)}{1 - \bar z_2 z_1}.
\end{equation}
Note that $k_{\mu, n}(z_1, z_2)$ is indeed an element of $\mathcal P_n$, i.e., a polynomial with respect to $z_1$ of degree at most $n-1$. It will be convenient to use a different representation of the reproducing kernel. Take $n \ge 1$, $a \in \D$, and define 
\begin{equation}\label{eq4}
\phit_{n} = \frac{z\phi_{n-1} - 
\bar a \phi^*_{n-1}}{\sqrt{1-|a|^2}}, 
\qquad
\phit_{n}^* = \frac{\phi^{*}_{n-1} - z a \phi_{n-1}}{\sqrt{1-|a|^2}}.
\end{equation}
\begin{Lem}\label{l4}
For all $z_1, z_2 \in \C$, we have
\begin{equation}\label{eq3}
\ov{\phit_{n}^*(z_2)}\phit_{n}^*(z_1)
-
\ov{\phit_{n}(z_2)}\phit_{n}(z_1) = 
\ov{\phi_{n}^*(z_2)}\phi_{n}^*(z_1)
-
\ov{\phi_{n}(z_2)}\phi_{n}(z_1)
\end{equation}
\end{Lem}
\beginpf The proof is a direct computation. At first, let $z_1 = z_2 = z$. Then the left hand side in \eqref{eq3} is equal to 
\begin{align*}
|\phit_{n}^*(z)|^2 -
|\phit_{n}(z)|^2
&=
\frac{1}{1-|a|^2}\left(|\phi_{n-1}^*(z) - za\phi_{n-1}(z)|^2 - 
|z\phi_{n-1}(z) - \bar a\phi_{n-1}^*(z)|^2\right) \\
&= 
|\phi_{n-1}^*(z)|^2 - 
|z\phi_{n-1}(z)|^2,
\end{align*}
which does not depend on $a$. Taking $a$ to be the recurrence coefficient $a_{n-1}$, we see that
\begin{equation}\label{eq6}
|\phit_{n}^*(z)|^2 -
|\phit_{n}(z)|^2 = |\phi_{n}^*(z)|^2 - 
|\phi_{n}(z)|^2.
\end{equation}
This relation holds for all $z \in \C$. Since functions in \eqref{eq3} are analytic in $z_1$ and anti-analytic in $z_2$, the lemma follows. \qed

\medskip

The following lemma is Corollary 4~in \cite{BD21}.
\begin{Lem}\label{l5}
For every $\lambda \in \D$ there is $a \in \D$ such that the corresponding polynomial $\phit_{n}^*$ in \eqref{eq4} defines a probability measure $\nu_{n,\lambda} = |\phit_{n}^*|^{-2}\,dm$ on $\T$ such that
\begin{equation}\label{eq5}
\K_{\nu_{n,\lambda}}(\lambda) \le \K_{\mu}(\lambda).
\end{equation}
%
\end{Lem}
In the rest of the paper, we use notation $\phit_{n}^*$ for the polynomial from Lemma \ref{l5}, where the value of the parameter $\lambda \in \D$ will be clear from the context.
\begin{Lem}\label{l6}
Let $\lambda \in \D$, and let $\phit_{n}^*$ be the corresponding polynomial from Lemma~\ref{l5}. We have 
\begin{equation}\label{eq17}
\int_{\T}\left|\frac{\phit_{n}^*(\lambda)}{\phit_{n}^*(\xi)} - 1\right|^2 \frac{1-|\lambda|^2}{|1 - \bar\xi \lambda|^2}\,dm(\xi) = e^{\K_{\nu_{n,\lambda}}(\lambda)} - 1.
\end{equation}
\end{Lem}
\beginpf By \eqref{eq12} and Lemma \ref{l4}, we have
\begin{equation}\label{eq7}
\frac{|\phit_{n}^*(z)|^2 -
|\phit_{n}(z)|^2}{1-|z|^2} = \sum_{k=0}^{n-1}|\phi_{k}(z)|^2,  \qquad z \in \C.
\end{equation}
It follows that the function $|\phit_{n}^*|^2 -
|\phit_{n}|^2$ is positive in $\D$ and is comparable to $1-|z|$ when $z$ approaches $\T$. Therefore, $\phit_{n}^*$ has no zeroes in $\D$. In fact, it has no zeroes in $\ov{\D} = \D \cup \T$ (if $\phit_{n}^*(z_0) = 0$ at some $z_0 \in \T$, then $1-|z| \lesssim |\phit_{n}^*|^2 - |\phit_{n}|^2 \lesssim |\phit_{n}^*|^2 \lesssim |z - z_0|^2$ near $z_0$, leading to a contradiction). 
It follows that the function $z \mapsto \frac{\phit_n^*(\lambda)}{\phit_n^*(z)}$ is analytic in a neighbourhood of $\ov{\D}$. Then the Poisson formula 
\begin{equation}\label{eq13}
u(\lambda) = \int_{\T}u(\xi) \frac{1-|\lambda|^2}{|1 - \bar\xi \lambda|^2}\,dm(\xi)
\end{equation}
for harmonic functions implies
\begin{equation}\label{eq14}
\int_{\T}\left|\frac{\phit_{n}^*(\lambda)}{\phit_{n}^*(\xi)} - 1\right|^2 \frac{1-|\lambda|^2}{|1 - \bar\xi \lambda|^2}\,dm = 
\int_{\T}\left|\frac{\phit_{n}^*(\lambda)}{\phit_{n}^*(\xi)}\right|^2 \frac{1-|\lambda|^2}{|1 - \bar\xi \lambda|^2}\,dm - 1,
\end{equation}
after noting that the function $u = -2\Re\Bigl(\frac{\phit_{n}^*(\lambda)}{\phit_{n}^*}\Bigr) + 1$ is harmonic in a neighbourhood of~$\ov{\D}$, $u(\lambda) = -1$. Observe that
\begin{equation}\label{eq15}
\int_{\T}\left|\frac{\phit_{n}^*(\lambda)}{\phit_{n}^*(\xi)}\right|^2 \frac{1-|\lambda|^2}{|1 - \bar\xi \lambda|^2}\,dm = 
|\phit_{n}^*(\lambda)|^2 \int_{\T}\frac{1-|\lambda|^2}{|1 - \bar\xi \lambda|^2}\,d\nu_{n,\lambda} = e^{\K_{\nu_{n,\lambda}}(\lambda)},
\end{equation}
because $|\phit_{n}^*(\lambda)|^2 = \exp\bigl(\int_{\T}\log |\phit_{n}^*(\xi)|^2 \frac{1-|\lambda|^2}{|1 - \bar\xi \lambda|^2}\,dm\bigr)$ (we use again formula~\eqref{eq13}, this time -- for the harmonic function $u = \log |\phit_{n}^*|^2$). The lemma now follows from \eqref{eq14} and \eqref{eq15}. \qed

\medskip

\begin{figure}
  \centering 
  \includegraphics[width=\textwidth]{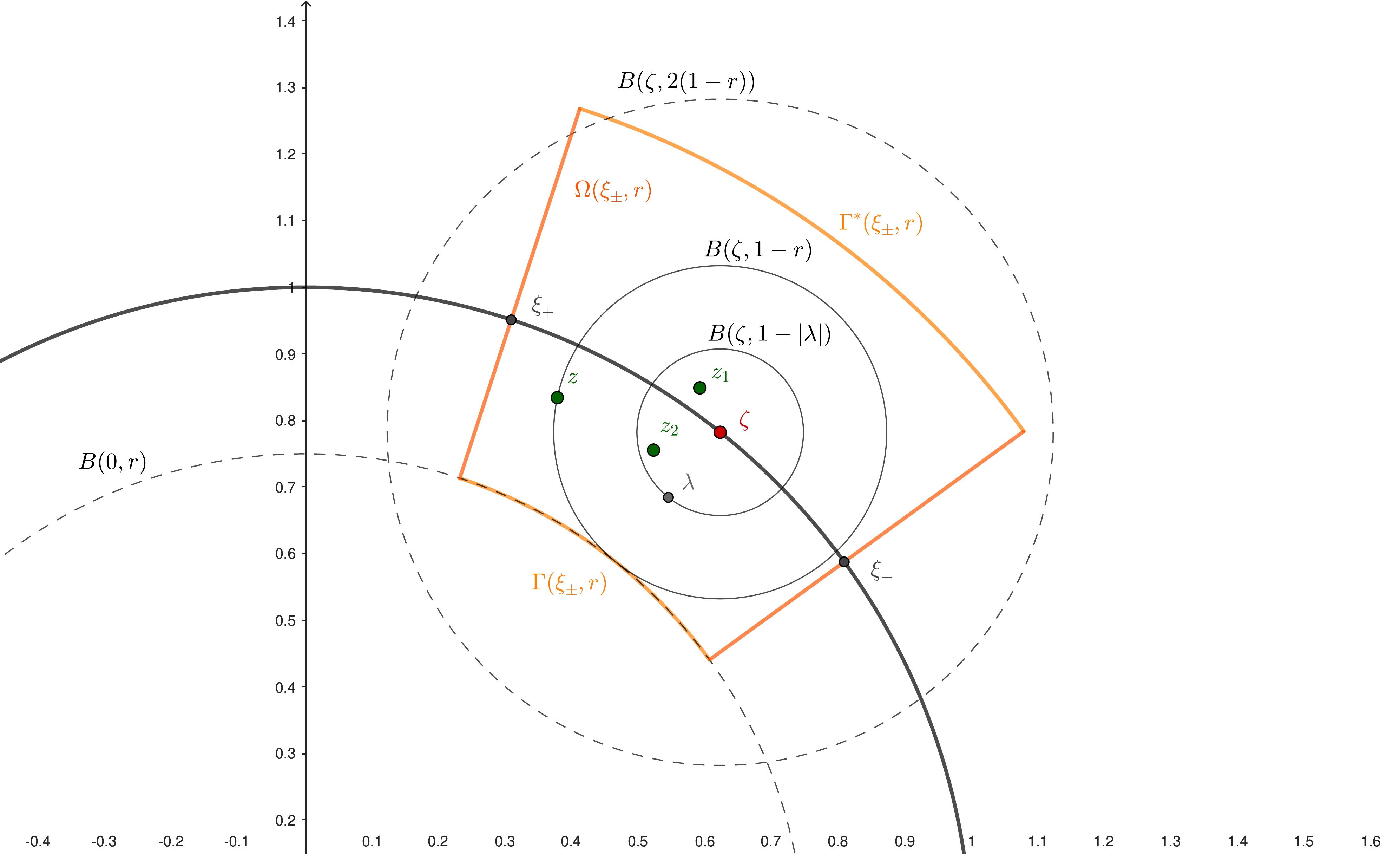}
  \caption{Objects that appear in the proof of Theorem \ref{t1}. Here $\xi_- = e^{i\pi/5}$, $\xi_+ = e^{2i\pi/5}$, $r =3/4$, $\zeta = e^{1.43\pi/5}$, $\lambda = (1+r)\zeta/2.$}  \label{fig1}
\end{figure}

Given two points $\xi_{\pm} \in \T$, $|\xi_+ - \xi_-| < 2$, and a number $r<1$, we denote by $\Gamma(\xi_{\pm}, r)$ the path in $\D$ formed by the union of two line segments 
$\{\rho \xi_{\pm}, \; \rho \in [r, 1)\}$ and the smaller arc of the circle $|z| = r$ with endpoints $r\xi_-$, $r\xi_{+}$. We also let $z^* = 1/\bar z$ for $z \in \C\setminus\{0\}$, and 
$$
\Gamma^*(\xi_{\pm}, r) = \{z \in \C:\; z^* \in \Gamma(\xi_{\pm}, r)\}.
$$
The union $\Gamma(\xi_{\pm}, r) \cup \{\xi_\pm\} \cup \Gamma^*(\xi_{\pm}, r)$ is then the boundary of a domain to be denoted by $\Omega(\xi_{\pm}, r)$. See Figure \ref{fig1}.

\begin{Lem}\label{l1}
Suppose that $h \in L^1(m)$, $\eta >0$, and $\lambda \in \D\setminus B(0, 3/4)$ are such that 
$$
\int_{\T}h(\xi) \frac{1-|\lambda|^2}{|1 - \bar\xi \lambda|^2}\,dm(\xi)
 \le \eta.
$$
Then there are $\xi_{\pm} \in \T$ such that for every $f \in H^1$ satisfying $|f| \le h$ on $\T$ we have $|f(z)| \lesssim \eta$ on $\Gamma(\xi_{\pm}, r)$, where $r$: $1-r = 2(1-|\lambda|)$. Moreover, $\xi_{\pm}$ are such that $\zeta = \lambda/|\lambda|$ belongs to the arc of $\T$ with endpoints $\xi_{\pm}$, and $1-r \le  |\xi_{\pm}-\zeta| \le 2(1-r)$.   
\end{Lem}
\beginpf Let $\lambda \in \D\setminus B(0, 3/4)$, $r=2|\lambda|-1$, $\zeta = \lambda/|\lambda|$. Consider the arc of the unit circle
$G = \T \cap B(\zeta,2(1-r))$, and define $m_{G} = \frac{\chi_{G}}{m(G)}dm$. We have
$$
\int_{\T} h \,dm_G \lesssim \eta.
$$
The set $G$ is the metric space with respect to the usual distance in $\C$. The measure $m_G$ has doubling property on this space: $m_G(B(\xi, 2\rho)) \le 4 m_G(B(\xi, \rho))$ for every $\xi \in G$, $\rho > 0$. It follows that the weak norm of Hardy-Littlewood maximal operator on $L^1(G, m_G)$ is bounded by a constant that does not depend on $G$. In other words, for every $t>0$, $g \in L^1(G, m_G)$ we have
$$
m_G\left(
\left\{\xi: \sup_{\rho>0}\frac{1}{m_G(B(\xi, \rho))}\int_{B(\xi, \rho)} g \,dm_G > t\right\}
\right)  \lesssim \frac{\|g\|_{L^1(G, m_G)}}{t},
$$
where the constant involved does not depend on $G$, $t$, $g$. Taking $g = h$, $t = \eps^{-1}\eta$ for some $\eps > 0$, we obtain
$$
m_G\left(
\left\{\xi: \sup_{\rho>0}\frac{1}{m_G(B(\xi, \rho))}\int_{B(\xi, \rho)} h \,dm_G > \eps^{-1}\eta \right\}
\right)  \lesssim \eps.
$$
It follows that
\begin{equation}\label{eq16}
m_G\left(
\left\{\xi: \sup_{\rho \in [0, 1)}\int_{G} h(u) \frac{1-\rho^2}{|1 - \rho\bar u \xi|^2}\,dm(u) > \eps^{-1}\eta\right\}
\right)  \lesssim \eps.
\end{equation}
Indeed, this follows from the fact that for each $\xi \in \T$, $\rho \ge 0$, the Poisson kernel $u \mapsto \frac{1-\rho^2}{|1 - \rho\bar u \xi|^2}$ can be uniformly approximated on $\T$ by positive convex combinations of functions of the form $\frac{\chi_{B(\xi, \delta)}}{m(B(\xi, \delta))}$, $\delta > 0$. Therefore, if $\xi \in G$ is such that for some $\rho \in [0, 1)$ we have
$$
\int_{G}h(u) \frac{1-\rho^2}{|1 - \rho\bar u \xi|^2}\,dm(u) > \eps^{-1}\eta
$$
then
$$
\frac{1}{m(B(\xi, \delta))}\int_{B(\xi, \delta)}\chi_{G}(u) h(u)\,dm(u) > \eps^{-1}\eta
$$
for some $\delta > 0$, and so  
$$
\frac{1}{m_G(B(\xi, \delta))}\int_{B(\xi, \delta)} h(u)\,dm_G(u) > \eps^{-1}\eta,
$$
proving \eqref{eq16}. Let us now take $\eps \in (0, 1)$ so small that the left hand side of \eqref{eq16} does not exceed $1/10$. Then there are $\xi_{\pm} \in G$ such that $\zeta$ belongs to the arc of $\T$ with endpoints $\xi_{\pm}$, we have $|\xi_{\pm} - \zeta| \ge 1-r$, and, moreover,
$$
\sup_{\rho \in [r, 1)}\int_{G} h \frac{1-\rho^2}{|1 - \rho\bar u \xi_{\pm}|^2}\,dm \le \eps^{-1}\eta.
$$   
For $u \in \T \setminus G$, we have 
$$
\sup_{\rho \in [r, 1)}\frac{1-\rho^2}{|1 - \rho\bar u \xi_{\pm}|^2}\,dm \lesssim \frac{1- |\lambda|^2}{|1 - \bar u \lambda|^2}.
$$
It follows that
\begin{equation}\label{eq19}
|f(\rho\xi_{\pm})| \lesssim \int_{\T\setminus G} h \frac{1- |\lambda|^2}{|1 - \bar u \lambda|^2}\,dm + \int_{G} h \frac{1-\rho^2}{|1 - \rho\bar\xi \xi_{\pm}|^2}\,dm \lesssim (1+\eps^{-1})\eta,
\end{equation}
for every $\rho \in [r, 1)$ with absolute constants. We also note that 
$$
\sup_{\xi \in G}\frac{1-r^2}{|1 - r\bar u \xi|^2} \lesssim \frac{1- |\lambda|^2}{|1 - \bar u \lambda|^2}, \qquad u \in \T,
$$
therefore, 
\begin{equation}\label{eq20}
|f(r\xi)| \lesssim \int_{\T} h(u) \frac{1- |\lambda|^2}{|1 - \bar u \lambda|^2}\,dm(u) \lesssim \eta,
\end{equation}
for $\xi \in G$. Collecting \eqref{eq19} and \eqref{eq20}, we see that $|f(z)| \lesssim \eta$ for $z \in \Gamma(\xi_{\pm}, r)$. \qed

\medskip

Define
$
\bt_n  = \frac{\phit_n}{\phit_n^*}. 
$
Formula \eqref{eq7} shows that $|\bt_n|\le 1$ on $\D$. In fact, $\bt_n$ is a Blaschke product of order $n$.
\begin{Lem}\label{l7}
Let $\lambda \in \D$, denote by $\eta = e^{\K_{\nu_{n,\lambda}}(\lambda)} - 1$ the number in the right hand side of \eqref{eq17}.  
Set $\alpha = \ov{\phit_n^*(\lambda)}/\phit_n^*(\lambda)$. We have
$$
\int_{\T}|\bt_n(\xi) - \alpha \xi^n|^2\frac{1-|\lambda|^2}{|1 - \bar \xi \lambda|^2}\,dm \lesssim \eta.
$$
\end{Lem}
\beginpf  Consider the sets 
$$
E_{l} = \left\{\xi \in \T: \; 
\left|\frac{\tilde\phi_n^*(\lambda)}{\tilde\phi_n^*(\xi)} - 1\right|^2 > 1/4\right\}, 
\qquad E_s = \T\setminus E_l.
$$ 
On $E_l$, the difference $|\bt_n(\xi) - \alpha \xi^n| \le 2$ could be large, but the measure of this set is small. Let us use Chebyshev's inequality and Lemma \ref{l6} to  estimate the corresponding integral:
$$
\int_{E_l}|\bt_n(\xi) - \alpha \xi^n|^2\frac{1-|\lambda|^2}{|1 - \bar \xi \lambda|^2} \le 
4\int_{E_l}\frac{1-|\lambda|^2}{|1 - \bar \xi \lambda|^2}\,dm \lesssim \eta.
$$
On $E_s$, the difference $|\bt_n(\xi) - \alpha \xi^n|$ is small. Indeed, for $\xi \in E_s$, we write
$$
|\bt_n - \alpha \xi^n| 
= \left|\frac{\bar\xi^n\phit_n(\xi)}{\phit_n^*(\xi)} - \alpha
\right|
=
\left|\ov{\frac{\phit_n^*(\xi)}{\ov{\phit_n^*(\lambda)}}} \frac{\phit_n^*(\lambda)}{\phit_n^*(\xi)} - 1
\right|.
$$
For $z \in \C$ such that $|1-z| < 1/2$, we have $|z^{-1}| \le 2$, so $|\phit_n^*(\xi)/\phit_n^*(\lambda)| \le 2$ on $E_s$. 
Then
$$
\left|\frac{\phit_n^*(\xi)}{\phit_n^*(\lambda)} - 1\right| = \left|\frac{\phit_n^*(\xi)}{\phit_n^*(\lambda)}\right|  \left|\frac{\phit_n^*(\lambda)}{\phit_n^*(\xi)} - 1\right| \le 2 \left|\frac{\phit_n^*(\lambda)}{\phit_n^*(\xi)} - 1\right|.
$$
Using $ab -1 = a(b-1) + a-1$ for $a = \ov{\phi_n^*(\xi)/\phi_n^*(\lambda)}$, $b = \phi_n^*(\lambda)/\phi_n^*(\xi)$, we see that
$$
|\bt_n - \alpha \xi^n| \lesssim \left|\frac{\tilde\phi_n^*(\lambda)}{\tilde\phi_n^*(\xi)} - 1\right|
$$
on $E_s$. The claim now follow from Lemma \ref{l6}. \qed
\begin{Lem}\label{l18}
Let $\lambda \in \D$, $A \ge 1$, $\eta = e^{\K_{\nu_{n,\lambda}}(\lambda)} - 1$. Define $r$ so that $1-r = 2(1-|\lambda|)$, and assume that $1/2 \le 1- A/n \le r < 1$.  There exists a number $\eta_0 \in (0,1)$ depending only on $A$, such that if $\eta \le \eta_0$, then there are  $\xi_{\pm} \in \T$ such that 
\begin{align}
|\bt_n(z) - \alpha z^n|
&\lesssim e^{4A}\sqrt{\eta}, \label{eq1}\\
\left|\frac{\phit_n^*(\lambda)}{\phit_n^*(z)} -1 \right| &\lesssim e^{2A}\sqrt{\eta}, \label{eq2}
\end{align}
for all $z \in \Gamma(\xi_{\pm}, r) \cup \Gamma^*(\xi_{\pm}, r)$. Moreover, $\xi_{\pm}$ are such that $\zeta = \lambda/|\lambda|$ belongs to the arc of~$\T$ with endpoints $\xi_{\pm}$, and $1-r \le |\xi_{\pm}-\zeta| \le 2(1-r)$.  
\end{Lem}
\beginpf Consider the function $h = |\bt_n - \alpha \xi^n|^2 + \left|\frac{\phit_n^*(\lambda)}{\phit_n^*} -1 \right|^2$ on $\T$. 
By Lemma~\ref{l7} and Lemma~\ref{l6}, we have
$$
\int_{\T}h(\xi)\frac{1-|\lambda|^2}{|1 - \bar \xi \lambda|^2}dm(\xi) \lesssim \eta.
$$
Then, by Lemma \ref{l1}, applied to the functions $(\bt_n(z) - \alpha z^n)^2$, $(\tfrac{\phit_n^*(\lambda)}{\phit_n^*(z)} -1)^2$,  there exists a contour $\Gamma(\xi_{\pm}, r)$ such that
\begin{equation*}
|\bt_n(z) - \alpha z^n|^2 
\lesssim \eta, \qquad 
\left|\frac{\phit_n^*(\lambda)}{\phit_n^*(z)} -1 \right|^2 \lesssim \eta,
\end{equation*}
for all $z \in \Gamma(\xi_{\pm}, r)$. Moreover, $\xi_{\pm}$ are such that $\zeta$ belongs to the arc of~$\T$ with endpoints $\xi_{\pm}$, and $1-r \le |\xi_{\pm}-\zeta| \le 2(1-r)$.  Choosing $\eta_0 \in (0,1)$ sufficiently small, one can guarantee that the left hand sides in these inequalities are smaller than $e^{-4A}/4$ for all $z \in \Gamma(\xi_{\pm}, r)$. In particular, for all $z \in \Gamma(\xi_{\pm}, r)$ we have
\begin{align*}
|z^n| &\ge r^n \ge (1-A/n)^n \ge e^{-2A},\\
|\bt_n(z)| &\ge|z^n| - |\alpha z^n - \tilde b_n| \ge e^{-2A} - e^{-2A}/2 \ge e^{-2A}/2,
\end{align*}
where we have used the elementary inequality $\log(1-x) \ge -2x$, $x \in [0, 1/2]$.
Then the identity
$$
\theta(1/\bar z) = \ov{1/\theta(z)}, \qquad z \in \C,
$$ 
for the inner functions $\theta = \bt_n$, $\theta = \alpha z^n$, gives \eqref{eq1} on $\Gamma^*(\xi_{\pm}, r)$:
$$
|\bt_n(1/\bar z) - \alpha (1/\bar z)^{n}| =
\frac{|\bt_n(z) - \alpha z^{n}|}{|\bt_n(z)z^{n}|} \lesssim e^{4A}\sqrt{\eta}, \qquad z \in \Gamma(\xi_{\pm}, r).
$$
To estimate $\bigl|\frac{\phit_n^*(\lambda)}{\phit_n^*} - 1\bigr|$ on $\Gamma^*(\xi_{\pm}, r)$, we use relation
$$
\frac{\phit_n^*(\lambda)}{\phit_n^*(1/\ov{z})} 
= \frac{\phit_n^*(\lambda)}{\ov{\phit_n(z)/z^n}}
= \frac{\bar z^n}{\ov{\bt_n(z)}}\frac{\ov{\phit_n^*(\lambda)}}{\ov{\phit_n^*(z)}}\frac{\phit_n^*(\lambda)}{\ov{\phit_n^*(\lambda)}} = 
\ov{\left(\frac{\phit_n^*(\lambda)}{\phit_n^*(z)}
\frac{\alpha z^n}{\bt_n(z)}\right)}.
$$
Then, formula $\ov{ab} - 1 = \ov{a(b-1) + a-1}$ for $a = \frac{\phit_n^*(\lambda)}{\phit_n^*(z)}$, $b=\frac{\alpha z^n}{\bt_n(z)}$ implies \eqref{eq2} on $\Gamma^*(\xi_{\pm}, r)$. \qed

\medskip

\begin{Lem}\label{l9}
Let $\lambda, A, n, r, \zeta, \eta$ be as in Lemma \ref{l18} and let $z_{1,2} \in B(\zeta, 1-|\lambda|)$. Then
$$
\left|\ov{\bt_n(z_1)}\bt_n(z_2) - \ov{z_1}^n z_2^n\right| \lesssim e^{4A}\sqrt{\eta} \cdot |1 - \bar z_1 z_2| \cdot n.
$$
\end{Lem}
\beginpf By Lemma \ref{l18}, relations \eqref{eq1}, \eqref{eq2} hold on some contour $\Gamma(\xi_{\pm}, r) \cup \Gamma^*(\xi_{\pm}, r)$, where $\xi_{\pm}$ are such that $\zeta = \lambda/|\lambda|$ belongs to the arc of~$\T$ with endpoints $\xi_{\pm}$, and $1-r \le |\xi_{\pm}-\zeta| \le 2(1-r)$. By the maximum modulus principle, \eqref{eq1} holds in the bounded domain $\Omega(\xi_{\pm}, r)$ such that 
$\partial \Omega(\xi_{\pm}, r) = \Gamma(\xi_{\pm}, r) \cup \Gamma^*(\xi_{\pm}, r) \cup \{\xi_{\pm}\}$. In particular, \eqref{eq1} holds for all $z \in B(\zeta, 1-r)$.  Now pick two points $z_1, z$ such that
$$
|z_1-\zeta| \le 1-|\lambda|, \qquad |z-\zeta| = 2(1-|\lambda|) = 1-r.
$$ 
We have
\begin{align*}
\left|\frac{\ov{\bt_n(z_1)}\bt_n(z) - \ov{z_1}^n z^n}{1-\bar z_1 z}\right| 
&\lesssim
\frac{\bigl|
(\bt_n(z_1) - \alpha z^n_1)(\bt_n(z)- \alpha z^n)\bigr|}{1-r} + \\
&+ \frac{\bigl|\ov{z_1}^n(\bt_n(z)- \alpha z^n)\bigr|}{1-r} + \frac{\bigl|(\alpha z_1^n - \bt_n(z_1)) z^n
\bigr|}{1-r}.
\end{align*}
Note that $|1-\bar z_1 z| \gtrsim 1-|\lambda| \ge \frac{1-r}{2} \gtrsim n^{-1},$ because $A \ge 1$. Recall that the maximum principle and \eqref{eq1} imply
$$
|\bt_n(z_1) - \alpha z^n_1| \lesssim e^{2A}\sqrt{\eta}, 
\qquad |\bt_n(z) - \alpha z^n| \lesssim e^{2A}\sqrt{\eta}.
$$
Next, we have
$$
\max(|z_1|^n,|z|^n) \le (1 + 1-r)^n \le (1 + A/n)^{n} \le e^{A}.
$$
Since $\eta \le \eta_0 \le 1$, we can conclude that
$$
\left|\frac{\ov{\bt_n(z_1)}\bt_n(z) - \ov{z_1}^n z^n}{1-\bar z_1 z}\right|  
\lesssim e^{4A}\sqrt{\eta} \cdot n.
$$
The function $z \mapsto \frac{\ov{\bt_n(z_1)}\bt_n - \ov{z_1}^n z^n}{1-\bar z_1 z}$ is analytic, hence the same estimate holds for all $z \in B(\zeta, 1-r)$ by the maximum modulus principle. In particular, it holds for all points $z_{1}$, $z_2$ in $B(\zeta, 1-|\lambda|)$. The lemma follows. \qed

\medskip

\noindent {\bf Proof of Theorem \ref{t1}.} Let $\mu \in \szc$, $\zeta \in \T$, $A \ge 1$. 
Let us choose $\tilde N_{0}$ such that $A/\tilde N_{0} < 1/4$ and
\begin{equation}\label{eq21}
\sup_{\rho \in [1 - A/\tilde N_{0}, 1)}e^{\K_{\mu}(\rho\zeta)} - 1 < \eta_0,
\end{equation}
where $\eta_0 \in (0, 1)$ is the number in Lemma \ref{l18}. For $n \ge \tilde N_{0}$, consider $z_{1,2} \in B(\zeta, A/n)$ and let $\lambda \in \D$ be defined by
$$
\lambda = (1-\max\limits_{k=1,2}|z_{k} - \zeta|)\zeta.
$$
Define $\eta = e^{\K_{\nu_{n,\lambda}}(\lambda)} - 1$ and observe that
$\eta < \eta_0$ by \eqref{eq21} and Lemma \ref{l5}.
We are in assumptions of Lemma~\ref{l18}. Estimate \eqref{eq2} and the maximum modulus principle give
$$
\frac{\phit_n^*(z)}{\phit_n^*(\lambda)} = 
\frac{1}{1 + R_1(z)}, \qquad |R_1(z)| \lesssim e^{2A}\sqrt{\eta},
$$
for all $z \in B(\zeta, 1-|\lambda|)$.
By Lemma \ref{l9}, for all $z_{1}, z_2 \in B(\zeta, 1-|\lambda|)$ we have
\begin{align*}
\frac{1 - \ov{\bt_n(z_1)}\bt_n(z_2)}{1 - \ov{z_1} z_2}
&=\frac{1 - \ov{z_1}^n z_2^n}{1 - \ov{z_1} z_2} + R_2(z_1, z_2) \cdot n,
\end{align*}
with $|R_2(z_1, z_2)| \lesssim e^{4A}\sqrt{\eta}$, so
\begin{align*}
k_{\mu, n}(z_1, z_2) 
&= \left(\frac{1 - \ov{z_1^n} z_2^n}{n(1 - \ov{z_1} z_2)} +  R_2(z_1, z_2)\right) \cdot \frac{|\phit_{n}^*(\lambda)|^2 \cdot n}{(1+R_1(z_1))\ov{(1+R_1(z_2))}}.\\
\end{align*}
Taking $z_1 = z_2 = \zeta$, we get
$$
k_{\mu, n}(\zeta, \zeta) = \frac{1 +  R_2(\zeta,\zeta)}{|1+R_1(\zeta)|^2} \cdot |\phit_{n}^*(\lambda)|^2 \cdot n.
$$
It remains to write 
$$
\frac{k_{\mu, n}(z_1, z_2)}{k_{\mu, n}(\zeta, \zeta)}
= 
\frac{\frac{1 - \ov{z_1}^n z_2^n}{n(1 - \ov{z_1} z_2)} + R_2(z_1, z_2)}{1+R_2(\zeta,\zeta)}\frac{|1+R_1(\zeta)|^2}{(1+R_1(z_1))\ov{(1+R_1(z_2))}}.
$$
Note that 
$$
\max_{z \in B(\zeta, 1-|\lambda|)}|R_{1}(z)| + \max_{z_{1,2} \in B(\zeta, 1-|\lambda|)}|R_2(z_1, z_2)| \le ce^{4A}\sqrt{\eta},
$$
for an absolute constant $c$. We see that one can take $\eps_0\in (0, \eta_0)$ so that $ce^{4A}\sqrt{\eps_0} \le 1/10$, then the required estimate will hold.

\medskip

\section{Proof of Theorem \ref{t2}}\label{s3}

\medskip

\noindent{\bf Proof of Theorem \ref{t2}.} Let $I$ be an arc of $\T$ with center $\zeta$ such that $|w(\xi) - w(\zeta)| \sim |\xi - \zeta|^{s}$ for all $\xi \in I$ and some $s > 0$. By the assumption of the theorem, we have $w(\zeta) > 0$. Replacing $I$ by a smaller interval, we can assume that  
$$
\delta < w(\xi) < 2\delta,
\qquad 
|\log w(\xi) - \log w(\zeta)| \sim |\xi - \zeta|^{s},
$$ 
for all $\xi \in I$ and some $\delta > 0$. Denote $u = \log w - \log w(\zeta)$.  We have 
$$
\int_{I}e^{u(\xi)}\frac{1 - \rho^2}{|1 - \rho \bar \xi \zeta|^2}\,dm = \int_{I}\bigl(1 + u(\xi) + O(u^2)\bigr)\frac{1 - \rho^2}{|1 - \rho \bar \xi \zeta|^2}\,dm.
$$
We also have 
\begin{align*}
\K_{\mu}(\rho\zeta) 
&= \log\left(\int_{\T}e^{u(\xi)}\frac{1 - \rho^2}{|1 - \rho \bar \xi \zeta|^2}\,dm\right) - 
\int_{\T}u(\xi)\frac{1 - \rho^2}{|1 - \rho \bar \xi \zeta|^2}\,dm\\
&\le \int_{\T}(e^{u(\xi)}-1)\frac{1 - \rho^2}{|1 - \rho \bar \xi \zeta|^2}\,dm - 
\int_{\T}u(\xi)\frac{1 - \rho^2}{|1 - \rho \bar \xi \zeta|^2}\,dm\\
&\lesssim (1 - \rho) \int_{\T\setminus I}(|e^{u(\xi)}-1| + |u(\xi)|)\,dm + \int_{I}u^2(\xi)\frac{1 - \rho^2}{|1 - \rho \bar \xi \zeta|^2}\,dm.
\end{align*}
By construction, we have
\begin{equation}\label{eq25}
\int_{I}u^2(\xi)\frac{1 - \rho^2}{|1 - \rho \bar \xi \zeta|^2}\,dm \lesssim \int_{I}|\xi-\zeta|^{2s}\frac{1 - \rho^2}{|1 - \rho \bar \xi \zeta|^2}\,dm.
\end{equation}
Let us set $\eps = 1-\rho$, assume that $\eps< m(I)/2$, and estimate
$$
\int_{I \cap B(\zeta,\eps)}|\xi-\zeta|^{2s}\frac{1 - \rho^2}{|1 - \rho \bar \xi \zeta|^2}\,dm
\lesssim
\frac{1}{\eps}\int_{I \cap B(\zeta,\eps)}|\xi-\zeta|^{2s}\,dm \lesssim \eps^{2s},
$$
and
\begin{equation}\notag
\int_{I \setminus B(\zeta,\eps)}|\xi-\zeta|^{2s}\frac{1 - \rho^2}{|1 - \rho \bar \xi \zeta|^2}\,dm
\sim
(1-\rho)\int_{I \setminus B(\zeta,\eps)}|\xi-\zeta|^{2s-2}\,dm. 
\end{equation}
We have 
\begin{equation}\label{eq26}
\int_{I \setminus B(\zeta,\eps)}|\xi-\zeta|^{2s-2}\,dm 
\sim
\begin{cases}
1, & s \in (1/2, +\infty),\\
\log \frac{1}{\eps},  & s =1/2,\\
\eps^{2s-1}, & s \in (0, 1/2),
\end{cases}
\end{equation}
as $\eps \to 0$. The constants in these relations depend on $s$. We see that 
$$
\K_{\mu}(\rho\zeta) \lesssim 
\begin{cases}
1-\rho, & s \in (1/2, +\infty),\\
(1-\rho)|\log(1-\rho)|,  & s =1/2,\\
(1-\rho)^{2s}, & s \in (0, 1/2),
\end{cases}
$$ 
as $\rho \to 1$. 
On the other hand, for $\eps = 1-\rho$ we have
\begin{equation}\label{eq18}
\K_{\mu}(\rho\zeta) 
\gtrsim \int_{\T}|u(\xi) - c|^2\frac{1 - \rho^2}{|1 - \rho \bar \xi \zeta|^2}\,dm
\end{equation}
Let $c_1$, $c_2$ be such that $c_1|\xi - \zeta|^{s} \le u(\xi) \le c_2|\xi - \zeta|^{s}$ for all $\xi \in I$. Choose $\delta \in (0, 1)$ so small that 
$$
c_1(1-\delta)^s - c_2\delta^s > 0.
$$
From \eqref{eq18} we see that there are $\xi_1 \in B(\zeta, \delta\eps)$, $\xi_2 \in B(\zeta, \eps) \setminus B(\zeta, (1-\delta)\eps)$ such that 
$$
|u(\xi) - c|^2 \lesssim \K_{\mu}(\rho\zeta), \qquad \xi = \xi_{1, 2}.
$$
Then 
$$
|u(\xi_1) - u(\xi_2)| \lesssim \sqrt{\K_{\mu}(\rho\zeta)},
$$
and, simultaneously,
$$|u(\xi_1) - u(\xi_2)| \ge c_1((1-\delta)\eps)^{s} - c_2(\delta\eps)^{s} = \eps^{s}(c_1(1-\delta)^s - c_2\delta^s) \gtrsim \eps^{s}.
$$
We see that $\K_{\mu}(\rho\zeta) \gtrsim (1-\rho)^{2s}$ for all $\rho$ close enough to $1$. A more accurate estimate is needed for $s \ge 1/2$. Let $\delta > 0$ be such that 
$$
c_1(1+\delta)^{s}/2 - c_2 > 0.
$$
Consider $j \in \Z$ satisfying $(1+\delta)^{j+3} < |I|/2$, and let $\xi_1$, $\xi_{2}$ be such that
$$(1+\delta)^{j} < |\xi_1-\zeta| < (1+\delta)^{j+1}, 
\qquad (1+\delta)^{j+2} < |\xi_2-\zeta| < (1+\delta)^{j+3}.
$$
We claim that either $|u(\xi_2) - c| \gtrsim |\xi_2 - \zeta|^s$ or $|u(\xi_1) - c| \gtrsim |\xi_1 - \zeta|^s$ for all such $\xi_1$, $\xi_2$, $j$. Indeed, if $|u(\xi_2) - c| \le c_1|\xi_2 - \zeta|^s/2$, then
\begin{align*}
|u(\xi_1) - c| 
&\ge |u(\xi_2)| - |u(\xi_1)| - |u(\xi_2) - c|\\ &\ge c_1|\xi_2 - \zeta|^s - c_2|\xi_1 - \zeta|^s - c_1|\xi_2 - \zeta|^s/2\\ 
&\gtrsim c_1(1+\delta)^{js+2s}/2 - c_2(1+\delta)^{js+s} \\
&\gtrsim (c_1(1+\delta)^{s}/2 - c_2)(1+\delta)^{js+s}\\
&\gtrsim (1+\delta)^{js+s} \gtrsim |\xi_1 - \zeta|^{s},
\end{align*}
proving the claim. We see that for each $j \in \Z$ such that $(1+\delta)^{j+3} < |I|/2$, we have $|u(\xi) - c| \gtrsim |\xi - \zeta|^{s}$ on a half of the set $\T \cap B(\zeta, (1+\delta)^{j+3}) \setminus B(\zeta, (1+\delta)^{j})$. It follows that 
$$
\K_{\mu}(\rho\zeta) \gtrsim \int_{I}|\xi - \zeta|^{2s}\frac{1 - \rho^2}{|1 - \rho \bar \xi \zeta|^2}\,dm,
$$
and estimate \eqref{eq26} shows that 
$$
\K_{\mu}(\rho\zeta) \gtrsim 
\begin{cases}
1-\rho, & s \in (1/2, +\infty),\\
(1-\rho)|\log(1-\rho)|,  & s =1/2,\\
(1-\rho)^{2s}, & s \in (0, 1/2).
\end{cases}
$$
This ends the proof of the theorem. \qed

\section{Example: the Poisson kernel}\label{s4}
In the following example, the asymptotic behaviour of ratios of reproducing kernels could be explicitly analysed. 
\begin{Ex}
Let $\lambda \in \D \setminus\{0\}$. Consider the probability measure $\mu = \frac{1 - |\lambda|^2}{|1 -  \lambda \xi|^2}\,dm$ on the unit circle $\T$. For $\zeta \in \T$, $u, v \in \R$, $z_1 = \zeta e^{u/n}$, $z_2 = \zeta e^{v/n}$, we have
$$
\frac{k_{\mu, n}(z_1, z_2)}{k_{\mu, n}(\zeta, \zeta)} = \frac{1 - \ov{z_2}^n z_1^n}{n(1 - \ov{z_2} z_1)} + \delta_n(u, v), \qquad n \to \infty,
$$
where $\sup\limits_{|u|,\,|v| \le 1}|\delta_n(u, v)|$ is comparable to $1/n$.
\end{Ex}
\beginpf We have (see Section 1.6 in \cite{Simonbook}) 
$$
\phi_{n}(z) = \frac{z^n - \bar\lambda z^{n-1}}{\sqrt{1-|\lambda|^2}}, 
\qquad
\phi^*_{n}(z) = \frac{1 - \lambda z}{\sqrt{1-|\lambda|^2}}.
$$
It follows that 
$$
k_{\mu, n}(z, z) = \sum_{0}^{n-1}|\phi_k(z)|^2 = 
\sum_{0}^{n-1}\frac{|z^{k}- \bar\lambda  z^{k-1}|^2}{1 - |\lambda|^2} = \frac{|z - \bar\lambda|^2}{|z|^2(1 - |\lambda|^2)}\frac{1 - |z|^{2n}}{1-|z|^2}.
$$
Then 
$$
\frac{k_{\mu, n}(z, z)}{k_{\mu, n}(\zeta, \zeta)} = \frac{|z - \bar\lambda|^2}{|z|^2|\zeta - \bar\lambda|^2}\frac{1 - |z|^{2n}}{n(1-|z|^2)},
$$
and we see that 
$$
\frac{k_{\mu, n}(z, z)}{k_{\mu, n}(\zeta, \zeta)} - 
\frac{k_{m, n}(z, z)}{k_{m, n}(\zeta, \zeta)}
= \frac{1 - |z|^{2n}}{n(1-|z|^2)}\frac{|z - \bar\lambda|^2 - |z|^2|\zeta - \bar\lambda|^2}{|z|^2|\zeta - \bar\lambda|^2}
$$
is comparable to
\begin{align*}
|z - \bar\lambda|^2 - |z|^2|\zeta - \bar\lambda|^2 
&= (1 - |z|^2)(|\lambda|^2 + 2\Re(\lambda(\zeta|z|^2-z)))\\ 
&\gtrsim (1 - |z|^2)|\lambda|^2 \gtrsim \frac{1}{n},
\end{align*}
if $\lambda \neq 0$, $|\zeta-z| \sim 1-|z|  \sim 1/n$, and $n$ is large. It follows that 
$$\sup\limits_{|u|,\,|v| \le 1}|\delta_n(u, v)| \gtrsim 1/n.
$$ 
Now take $z_1 = \zeta e^{u/n}$, $z_2 = \zeta e^{v/n}$, with $|u|,\,|v| \le 1$.  We have
\begin{align*}
k_{\mu, n}(z_1, z_2) 
&= \frac{1}{1-|\lambda|^2}\frac{(1 - \lambda z_1)\ov{(1 - \lambda z_2)} - z_1^n \bar z_2^n(1 - \bar\lambda/z_1)(1 - \lambda/\bar z_2)}{1 - z_1 \bar z_2 }\\
&= \frac{(1 - \lambda z_1)\ov{(1 - \lambda z_2)}}{1-|\lambda|^2}\frac{1 - z_1^n \bar z_2^n\theta_\lambda(z_1, z_2)}{1 - z_1 \bar z_2 },
\end{align*}
where 
$$
\theta_{\lambda}(z_1, z_2) = \frac{(1 - \bar\lambda/z_1)(1 - \lambda/\bar z_2)}{\ov{(1 - \lambda z_2)}(1 - \lambda z_1)} = 1 + O(1/n).
$$
In particular, we have
\begin{align*}
\frac{k_{\mu, n}(z_1, z_2)}{k_{\mu, n}(\zeta, \zeta)} 
&= \frac{(1 - \lambda z_1)\ov{(1 - \lambda z_2)}}{|1 - \lambda \zeta|^2}\frac{1 - z_1^n \bar z_2^n\theta_\lambda(z_1, z_2)}{n(1 - z_1 \bar z_2)}.
\end{align*}
Here 
$$
\frac{(1 - \lambda z_1)\ov{(1 - \lambda z_2)}}{|1 - \lambda \zeta|^2} = 1 + O(1/n),
$$
and
$$
\frac{1 - z_1^n \bar z_2^n\theta_\lambda(z_1, z_2)}{n(1 - z_1 \bar z_2)} = \frac{k_{m, n}(z_1, z_2)}{k_{m, n}(\zeta, \zeta)}  + z_1^n \bar z_2^n\frac{\theta_\lambda(z_1, z_2)-1}{n(1 - z_1 \bar z_2)}.
$$
Note that $\theta_\lambda(z_1, z_2) =1$ in the case $z_1 \bar z_2 = 1$. Considering $z_1$, $z_2$ such that $|1 - z_1 \bar z_2| \sim 1/n$ and using maximum modulus principle for analytic functions, we see that 
$$
\left|z_1^n \bar z_2^n\frac{\theta_\lambda(z_1, z_2)-1}{n(1 - z_1 \bar z_2)}\right| \lesssim \frac{1}{n}
$$
for all $z_1$, $z_2$ such that $|z_{1,2} -\zeta| \lesssim 1/n$. It follows that $\sup\limits_{|u|,\,|v| \le 1}|\delta_n(u, v)| \le 1/n$. \qed

\section{Sharpness of Theorem \ref{t1}.}\label{s5}
Sharpness of Theorem \ref{t1} is an open problem. As we have seen in the previous section, one can have
$$
\sup_{|z_{1,2} - 1| \le 1/n}\left|\frac{k_{\mu, n}(z_1, z_2)}{k_{\mu, n}(1, 1)} - \frac{1 - \ov{z_2}^n z_1^n}{n(1 - \ov{z_2} z_1)}\right| =  o\left(\sqrt{\K_{\mu}(1-1/n)}\right), 
$$
if the measure $\mu$ is very regular. Our aim  in this section is to present some numerical results for the measures of the form $w_s\,dm$, where $w_s(e^{i\theta}) = c_{s}e^{|\theta|^{s}}$, $\theta \in [-\pi, \pi]$, and the constant $c_s$ is chosen so that $\int_{\T} w_s \,dm = 1$. We will see that 
$$
\sup_{|z_{1,2} - 1| \le 1/n}\left|\frac{k_{w_s, n}(z_1, z_2)}{k_{w_s, n}(1, 1)} - \frac{1 - \ov{z_2}^n z_1^n}{n(1 - \ov{z_2} z_1)}\right| \gtrsim \sqrt{\K_{w_s}(1-1/n)}, 
$$
in the case $s=0.1$, $s = 0.2$, and $s = 0.4$ (similar results can be obtained for other values of $s \in (0, 1/2)$). In particular, the upper bound in Theorem \ref{t1} in these cases coincides with the lower bound up to a multiplicative factor comparable to $1$. We use the following MATLAB code to produce our examples (note that the weight $w_s$ and the orthogonal polynomials in the script are not normalized, because the normalization plays no role when we consider ratios of reproducing kernels).

\medskip

\begin{verbatim}
%% Clear window and variables
clc
clear vars 

%% Initialization

% Maximal degree
N=8000;

% Step (N/s must be positive integer)
step = 20;

% Holder index
s = 0.4;

% Initial difference
diffold = 0;

% Auxiliary arrays
MMNTS=zeros(N,1);              % array for moments
D = zeros(1,N/step);           % array for differences
alphaCand = zeros(1,N/step);   % array for Holder index candidates   
CalphaCand = zeros(1,N/step);  % array for C_alpha candidates                     

%% Computing the moments of the weight e^{|x|^s}, x in [-pi, pi]
for j=0:1:N-1
fun = @(x) cos(j*x).*exp(abs(x).^s);     % the weight is even
MMNTS(j+1) = integral(fun,-pi,pi);
end

%% Cycle for comparing the ratios of kernels for different n
for n=step:step:N
disp(['n = ', num2str(n)]);
% Define the reference point
xn = 1-1/n;                                     
% Compute the free reproducing kernel of order n-1 at (xn, xn)
RepKer0atxnxn =((xn).^(0:1:n-2))*((xn).^(0:1:n-2))'; 
% Compute the free reproducing kernel of order n-1 at (xn, xn)
RepKer0at11 = n-1; 
% Define the Toeplitz matrix 
% and the basis vector e_n = (0,0..., 0,1)'
j=0:1:n-1;
moments = double(MMNTS(1:n));
T = toeplitz(moments);
eo = double(1:n == n)';
% Use MATLAB function invToeplitz() by John P. Cunningham for fast
% invertion of the Toeplitz matrix T to compute coeffitients of OPs
% see https://github.com/aecker/gpfa.git
coefOP=((invToeplitz(T))*eo)';                  
coefRefOP = conj(flip(coefOP));
% Compute value of OP and reflected OP at 1
OPat1 = coefOP*ones(n,1);
RefOPat1 = coefRefOP*ones(n,1);
% Compute value of OP and reflected OP at xn
OPatxn = coefOP*((xn).^(0:1:n-1))';
RefOPatxn = coefRefOP*((xn).^(0:1:n-1))';
% Compute coeffitients of reproducing kernel at 1
% Numerator
RepKerNumat1 = coefRefOP*conj(RefOPat1) - coefOP*conj(OPat1);         
% Denumerator
RepKerDenat1 = [1, -1];                     
% Kernel                          
[RepKerat1, r1] = deconv(flip(RepKerNumat1), flip(RepKerDenat1));     
% Compute coeffitients of reproducing kernel at xn
% Numerator
RepKerNumatxn = coefRefOP*conj(RefOPatxn) - coefOP*conj(OPatxn);      
% Denumerator
RepKerDenatxn = [1, -conj(xn)];             
% Kernel                          
[RepKeratxn,rxn]=deconv(flip(RepKerNumatxn),flip(RepKerDenatxn)); 
% Compute the value of reproducing kernel at (1,1)
RepKerat11 = flip(RepKerat1)*ones(n-1,1);
% Compute the value of reproducing kernel at (xn,xn)
RepKeratxnxn = flip(RepKeratxn)*((xn).^(0:1:n-2))';

%% Compare the ratios
ratio0 = double(RepKer0atxnxn/RepKer0at11);
ratio1 = double(RepKeratxnxn/RepKerat11);
diff = ratio0-ratio1;
D(n/step) = abs(diff);
% Assuming power-type law "C_alpha n^{-alpha}", find the candidate 
% for alpha for the current value of n
alphaCand(n/step) = (n/step).*(diffold./diff - 1); 
% Assuming power-type law "C_alpha n^{-alpha}", find the candidate 
% for C_alpha for the current value of n
CalphaCand(n/step) = diff.*(n.^alphaCand(n/step));
diffold = diff;
end
\end{verbatim}

\medskip

\noindent This script produces the array \verb|D| consisting of differences 
$$
\verb|D(n/20)| = \left|\frac{k_{w_s, n}(x_n, x_n)}{k_{w_s, n}(1, 1)} - \frac{1 - |x_n|^{2n}}{n(1 - |x_n|^2)}\right|, \qquad n = 20, 40, 60, \ldots, 8000,
$$
for the value $s = 0.4$ of H\"older continuity index. One can plot the values of this array (see below) and observe that it behaves like $C_{\alpha}n^{-\alpha}$ for certain values  $\alpha>0$, $C_{\alpha}>0$. To find these values, we observe that
$$
\frac{C_{\alpha}(n-20)^{-\alpha}}{C_{\alpha}n^{-\alpha}} = \left(1 - \frac{20}{n}\right)^{-\alpha} = 1 + \frac{20\alpha}{n} + O(n^{-2}),
$$ 
which gives a reasonable recipe to compute the array \verb|alphaCand| of candidates for the power index $\alpha$: 
$$
\verb|alphaCand(n/20)| = \frac{n}{20}\left(\frac{D(n/20-1)}{D(n/20)} - 1\right), \qquad 
n=20, 40, \ldots, 8000.
$$
When printed, the last $7$ elements of this array (corresponding to values $n=7880,7900$, \ldots, $8000$) look as follows:
$$
\ldots
\quad
0.39361 \quad 0.39362 \quad 0.39363 \quad   0.39363 \quad 0.39364 \quad 0.39365 \quad 
0.39365
$$
which is pretty close to $s = 0.4$. So it seems plausible that 
$$
\left|\frac{k_{w_s, n}(x_n, x_n)}{k_{w_s, n}(1, 1)} - \frac{1 - |x_n|^{2n}}{n(1 - |x_n|^2)}\right| \ge \frac{C_{0.4}}{n^{0.4}} 
$$
for all $n$ large enough. We can even suggest a candidate for the constant $C_{0.4}$ from our numerical data. For this we use approximation 
$$
C_{0.4} \sim C_{\alpha} = N^{0.39365}D(N/20) = 0.03791..., \qquad N = 8000.
$$ 
The plots of the resulting functions
$$f_1(n) = \left|\frac{k_{w_s, n}(x_n, x_n)}{k_{w_s, n}(1, 1)} - \frac{1 - |x_n|^{2n}}{n(1 - |x_n|^2)}\right|, 
\qquad f_2(n) = C_{0.4}n^{-0.4},
$$ 
are given below. We also consider the cases $s = 0.1$ and $s = 0.2$.
\newpage
\begin{figure}
  \centering 
  \includegraphics[width=\textwidth]{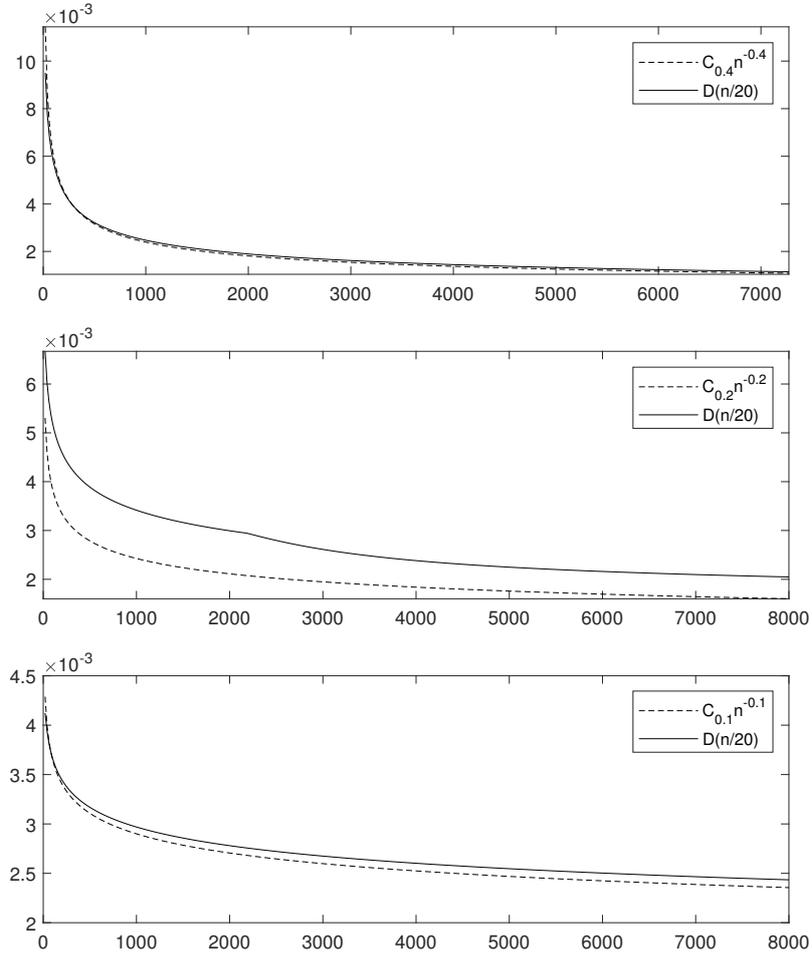}
  \caption{Graphs of functions $f_1$, $f_2$ for $s = 0.1$, $s = 0.2$, $s = 0.4$.}  
 \label{fig2}
\end{figure}
The graphs demonstrate the fact that $f_1 \ge f_2$ for $s \in \{0.1, 0.2, 0.4\}$ and large integers $n$. Similar numerical results can be obtained for other values $s \in (0,1/2)$.


\bibliographystyle{plain} 
\bibliography{bibfile}
\enddocument

\enddocument